\documentclass[11pt,reqno]{amsart}
\usepackage{style/paper}
\newcommand{\R}{\mathbb{R}}

\newcommand{\Q}{\mathbb{Q}}

\DeclareMathOperator{\supp}{supp}

\DeclareMathOperator{\ann}{ann}

\title[Inequalities from Hodge--Riemann relations]{Inequalities for intersection numbers from Hodge--Riemann relations}
\author{Xiangru Zeng}
\email{xiangru\_zeng@berkeley.edu}
\date{September 2026}
\begin{document}
\begin{abstract}
We derived inequalities for intersection numbers of nef classes inside algebras satisfying the K\"ahler package, using Hodge-Riemann relations up to degree two. As an application, we construct certain strictly Lorentzian polynomials and show that they cannot be realized as the volume polynomial of nef classes inside such an algebra.
\end{abstract}
\maketitle
\section{Introduction}\label{sec: introduction}

Let $A^\bullet=\bigoplus_{k=0}^dA^k$ be a finitely generated, commutative graded $\R$-algebra with a degree map $\deg\colon A^d\to \R$. $A^\bullet$ is said to satisfy the \textit{K\"ahler package} with respect to an open convex cone $K\subseteq A^1$ if
\begin{enumerate}
\item (Poincar\'e duality) $\deg$ is an isomorphism, and  for $0\leq k\leq d$, the bilinear pairing\[ A^{k}\times A^{d-k}\to\R\]that sends $(\alpha,\beta)$ to $\deg(\alpha\beta)$ is non-degenerate. In particular $A^0\cong \R$.
\item (mixed Hard Lefschetz) For $0\leq k\leq \lfloor\frac{d}{2}\rfloor$ and any $\ell_1,\ldots,\ell_{d-2k}\in K$, the linear map \[A^k\to A^{d-k},\alpha\mapsto \left(\prod_{i=1}^{d-2k}\ell_i\right)\alpha\]is an isomorphism.
\item (mixed Hodge-Riemann relations) For $0\leq k\leq \lfloor\frac{d}{2}\rfloor$ and any $\ell_0,\ell_1,\ldots,\ell_{d-2k}\in K$, the symmetric bilinear form \[A^k\times A^k\to \R,\quad (\alpha,\beta)\mapsto(-1)^k \deg\left(\left(\prod_{i=1}^{d-2k}\ell_i\right)\alpha\beta\right)\]is positive-definite on the kernel of the map\[A^k\to A^{d-k+1},\alpha\mapsto \left(\prod_{i=0}^{d-2k}\ell_i\right)\alpha.\]
\end{enumerate}
Classes inside $K$ are often called \textit{ample classes}, whereas classes inside $\overline K$ are called \textit{nef classes}. Prototypical examples of such an algebra would be the rings of real $(k,k)$-cohomology classes on a compact K\"ahler manifold.

The K\"ahler package imposes strong constraints on the intersection numbers of nef classes. The \textit{Alexandrov-Fenchel type} inequalities assert that,
for any classes $\alpha\in A^1$ and $\beta,\ell_3,\ldots,\ell_d\in \overline K$, we have \[
    \deg(\alpha\beta\ell_3\ldots\ell_d)^2\geq\deg(\alpha^2\ell_3\ldots\ell_d)\deg(\beta^2\ell_3\ldots\ell_d).
\]
These inequalities follow from Hodge-Riemann relations up to degree one.

By investigating Hodge-Riemann relations up to degree two, we discover a new family of inequalities. Let $R=\R[X_1,\ldots,X_n]$ equipped with the standard grading, and $K_R=\{\sum_{i=1}^n \lambda_iX_i\mid \lambda_i>0,\forall i\}$, we define partial orders on $R^1$ by saying $A\preceq B$ when $B-A\in \overline K_R$, and $A\prec B$ if $B-A\in K_R$. 
\begin{theorem}\label{thm: main}
Suppose that
\[
    HP=\sum_{i=1}^rA_iB_iV_i\in R\tag{1}\label{eq: referencing}
\]
where $P\in R^2$, $H, A_i,B_i,V_i\in R^1$ and $H,A_i,B_i\in \overline{K_R}$ and $0\preceq A_i,B_i\preceq H$. Let $A^\bullet$ be a graded algebra of degree $d\geq 4$ satisfying the K\"ahler package with respect to $K$. We consider a graded ring homomorphism $f\colon R\to A$ that maps $K_R$ to $\overline K$, and we denote the images of $P,H,A_i,B_i,V_i$ by $p,h,\alpha_i,\beta_i, v_i$ respectively.

For any class $\ell\in A^{d-4}$ given by the product of $(d-4)$ classes in $\overline{K}$, we have the following inequality
\[
\deg (p^2\ell)+25\left(\sum_{i=1}^r \sqrt{\frac{2\deg(\alpha_i\beta_i h v_i\ell)^2}{\deg(\alpha_i\beta_ih^2\ell)}-\deg(\alpha_i\beta_i v_i^2\ell)}\right)^2\geq \frac{2\deg(ph^2\ell)^2}{\deg(h^4\ell)}
\]
whenever $\deg(\alpha_i\beta_ih^2\ell)\neq 0$ for all $i$.
\end{theorem}

To show that these inequalities are not vacuous, we give an example showing that they do not follow from the AF-type inequalities.

\begin{example}\label{ex: ternary}
    Put
\[
        a=\frac{2023}{1500},\qquad
        b=\frac{119}{100},\qquad
        r=\frac{11}{10},\qquad
        s=\frac12,
\]
and define
\[
\begin{aligned}
F(x,y,z)={}&
 a\frac{x^4}{24}
+b\frac{x^3y+x^3z}{6}
+\frac{x^2y^2+x^2z^2+y^2z^2}{4}
+r\frac{x^2yz}{2}                                      \\
&+s\frac{xy^3+xz^3}{6}
+\frac{xy^2z+xyz^2}{2}
+s^2\frac{y^4+z^4}{24}
+s\frac{y^3z+yz^3}{6}.
\end{aligned}
\]
Then for any degree $4$ graded $\R$-algebra $A$ and classes $\chi,\gamma,\zeta\in A^1$ such that
\[
    \frac{1}{4!}\deg(x\chi+y\gamma+z\zeta)^4=F_\epsilon(x,y,z)\coloneqq F(x+\epsilon\sigma,y+\epsilon\sigma,z+\epsilon\sigma),\; \sigma=x+y+z,
\]
the inequalities in \Cref{thm: main} fail for the map $R=\R[X,Y,Z]\to A,\ g(X,Y,Z)\mapsto g(\chi,\gamma,\zeta)$ when $0<\epsilon\ll 1$, where the equality \eqref{eq: referencing} is given by
\begin{align*}
    H  =& X+Y+Z,\\
    P  =& 750X^2-750XY-750XZ-28Y^2-28Z^2+735YZ,\\
    HP =&
        \ X^2\cdot 25(30X-17Y-17Z)\\
        &+XY\cdot \frac{85}{2}(10X-2Y-9Z)+XZ\cdot \frac{85}{2}(10X-9Y-2Z)\tag {1'}\label{eq: referencing'}\\
        &-Y^2\cdot (693(X-Z)+14(2Y-Z))-Z^2\cdot(693(X-Y)+14(2Z-Y)).
\end{align*}
On the other hand, $F_\epsilon(x,y,z)$ is a \textit{strictly Lorentzian polynomial}, which means that it has positive coefficients and the Hessian matrices of its second order partial derivatives have the same signature as Lorentzian metrics. For these polynomials, one can canonically construct Poincar\'e duality algebras satisfying AF-type inequalities.
\end{example}

We organize this paper as follows. In section \ref{sec: norms}, we prove auxiliary results on norm estimates associated to symmetric forms with Lorentzian signature; in section \ref{sec: proof}, we will apply these results to an external inequality \eqref{eq: external inequality} to prove the main theorem; then in section \ref{sec: ex}, we will give detailed analysis on \cref{ex: ternary} and compare different sets of polynomials.

\subsection*{AI disclosure} The prototype of \cref{thm: main}, its proof, and \cref{ex: ternary} were discovered with the assistance of Danus, a multi-agent mathematical reasoning system \cite{LGS+26}. During the research process, Danus employed three OpenAI models: GPT-5.5, GPT-5.6 Sol, and GPT-6 Astra. Separately, GPT-6 Astra was used for literature search and proofreading the manuscript. The results shown in this paper were reshaped and strengthened by the author for better comprehensibility. The author takes full responsibility for the mathematical content and correctness of the paper.

\subsection*{Acknowledgments} The author thanks Christopher Eur for helpful comments. The author also thanks his advisor David Eisenbud for suggestions on writing.
\section{Estimates for Lorentzian forms}\label{sec: norms}

Let $A^\bullet$ be an algebra of degree $d\geq2$ and satisfy the
K\"ahler package with respect to $K$. For
$\ell_0,\ldots,\ell_{d-2}\in K$, define
\[
 T(x,y)=\deg(\ell_1\cdots\ell_{d-2}xy),\qquad x,y\in A^1.
\]
The degree-zero relations give $T(\ell_0,\ell_0)>0$, and the degree-one
relations say that $T$ is negative definite on the hyperplane
$\{x:T(\ell_0,x)=0\}$. Thus $T$ has signature
$(1,\dim A^1-1)$, where the first entry counts positive eigenvalues.
We call this the \textit{Lorentzian signature}.

Let $W$ be a finite-dimensional real vector space and $T$ be a symmetric bilinear form with the Lorentzian signature. For any $h\in W$ such that $T(h,h)>0$, we can associate a $T$-invariant involution
\[
    \theta_{T,h}\colon W\to W,\quad x\mapsto \frac{2T(h,x)}{T(h,h)}\cdot h-x,
\]
which fixes $h$ and negates $h^\perp$. We define
\[
    G_{T,h}(x,y)=T(\theta_{T,h}(x),y)=\frac{2T(h,x)T(h,y)}{T(h,h)}-T(x,y).
\]
Since $T$ is negative definite on $h^\perp$, for orthogonal decomposition $x=\lambda h+x^\perp$, we have
\[
G_{T,h}(x,x)=T(\lambda h-x^\perp,\lambda h+x^\perp)=\lambda^2T(h,h)-T(x^\perp,x^\perp)\geq 0,
\]
with equality only when $x=0$. Thus $G_{T,h}$ is an inner product on $W$.
\begin{lemma}\label{lem: cauchy-schwarz}
    \[
        T(x,y)^2\leq G_{T,h}(x,x)G_{T,h}(y,y),\quad \forall x,y\in W
    \]
\end{lemma}
\begin{proof}
    \begin{align*}
         T(x,y)^2=& G_{T,h}(\theta_{T,h}x,y)^2\\
         \leq& G_{T,h}(\theta_{T,h}(x),\theta_{T,h}(x))G_{T,h}(y,y)\\
         =&G_{T,h}(x,x)G_{T,h}(y,y),
    \end{align*}
    where the inequality is Cauchy--Schwarz for $G_{T,h}$.
\end{proof}
\begin{lemma}\label{lem: estimate of norms}
    Let $R,S,T$ be symmetric bilinear forms with Lorentzian signature on $W$ satisfying $T=R+S$. Suppose that $h\in W$ satisfies
    \[
        r=R(h,h)>0,\ s=S(h,h)>0.
    \]
    Further assume that there is $u\in W$ with
    \[
        R(h,x)=T(u,x),\quad S(h,x)=T(h-u,x)
    \]
    for all $x\in W$, and that $T(u,u),T(h-u,h-u)\geq 0$. Then we have
    \[
        G_{R,h}(x,x)\leq 5 G_{T,h}(x,x),\quad \forall x\in W.
    \]
\end{lemma}
\begin{proof}
    Set $t=r+s=R(h,h)+S(h,h)=T(h,h)$. We have
    \begin{align*}
        G_{T,h}(u,u)=\frac{2T(u,h)^2}{T(h,h)}-T(u,u)=\frac{2R(h,h)^2}{t}-T(u,u)\leq \frac{2 r^2}{t}.
    \end{align*}
    By \cref{lem: cauchy-schwarz}, we have
    \[
        R(h,x)^2=T(u,x)^2\leq G_{T,h}(u,u)G_{T,h}(x,x)\leq \frac{2r^2}{t}G_{T,h}(x,x).
    \]
    Similarly, we have
    \[
        S(h,x)^2\leq \frac{2s^2}{t}G_{T,h}(x,x) .
    \]
    Since $G_{S,h}$ is positive definite, we have
    \[
        S(x,x)=\frac{2S(h,x)^2}{s}-G_{S,h}(x,x)\leq \frac{2S(h,x)^2}{s}.
    \]
    Therefore
    \begin{align*}
        G_{R,h}(x,x)=&\frac{2R(h,x)^2}{r}-R(x,x)\\
        =&\frac{2R(h,x)^2}{r}+S(x,x)-T(x,x)\\
        \leq&\frac{4r+4s}{t}G_{T,h}(x,x)+\left(G_{T,h}(x,x)-\frac{2T(h,x)^2}{t}\right)\\
        \leq& 5G_{T,h}(x,x).
    \end{align*}
\end{proof}

\section{Proof of the main theorem} \label{sec: proof}
We first reduce to the case of degree $4$ algebras. Let $A^\bullet$ be a commutative graded $\R$-algebra equipped with a degree map $\deg_A\colon A^d\to\R$ satisfying the Poincaré duality.  For $x\in A$, let $\ann(x)$ denote the annihilator of $x$. When $x$ is nonzero and homogeneous of degree $k$, the quotient ring $A/\ann(x)$ is graded of degree $d-k$, with an induced degree map
\[
    \deg_{A/\ann(x)}(p(y))=\deg_A(xy)
\]
for all $y\in A^{d-k}$, where $p\colon A\to A/\ann(x)$ is the natural projection. We recall the following descent lemma.
\begin{proposition}[{\cite[Lemma 1.16]{CKS87}} and {\cite[Theorem 2.1.5]{KK87}}, see also {\cite[Theorem 3.2]{Cat08}} and {\cite[Lemma 4.5]{AHL25}}]\label{prop: descent}
    $A/\ann(x)$ satisfies Poincaré duality. If furthermore $A$ has the Kähler package with respect to an open convex cone $K_A$ and $x\in \overline K_A$, then $A/\ann(x)$ satisfies the Kähler package with respect to the image of $K_A$.
\end{proposition}

By repeatedly applying \cref{prop: descent}, we see that for $\ell_5,\ldots,\ell_d\in \overline{K}_A\setminus 0$, whose product is nonzero, we have
\[
    \deg_{A/\ann(\ell_5\cdots\ell_d)}(x)=\deg_{A}(x \ell_5\cdots\ell_d)
\]
and that $A/\ann(\ell_5\cdots\ell_d)$ satisfies the K\"ahler package with respect to the projection of $K_A$. To prove \cref{thm: main} for $A$ with respect to $K_A$, it suffices to prove it for $A/\ann(\ell_5\cdots\ell_d)$ with respect to the image of $K_A$.

For the rest of this section, we assume $A=\bigoplus_{i=0}^4 A^i$ which has the K\"ahler package with respect to degree map $\deg$ and open convex cone $K\subseteq A^1$. For any $\alpha,\beta\in K$, we have the Lorentzian bilinear form on $A^1$
\[
    T_{\alpha\beta}(x,y)\coloneqq \deg(\alpha\beta xy).
\]
For any $h\in K$, we can further define the inner product
\[
    G_{\alpha\beta,h}(x,y)\coloneqq G_{T_{\alpha\beta},h}(x,y)=\frac{2 \deg(\alpha\beta h x)\deg(\alpha\beta h y)}{\deg(\alpha\beta h^2)}-\deg(\alpha\beta x y).
\]
We partially order $A^1$ by saying $\alpha \prec \beta$  (or $\alpha \preceq \beta$) when $\beta-\alpha\in K$ (or $\beta-\alpha\in \overline{K}$). Then we have the following interpretation of \cref{lem: estimate of norms}.
\begin{lemma} \label{lem: estimate of norms revisited}
    Let $\alpha,\beta,\gamma\in K$ such that $\alpha\preceq\beta$. Then for every $x\in A^1$,
    \[
        G_{\alpha \gamma,\beta}(x,x)\leq 5G_{\beta \gamma,\beta}(x,x).
    \]
\end{lemma}
\begin{proof}
    In the case where $\alpha\prec \beta$, let $R(x,y)=T_{\alpha \gamma}(x,y), S(x,y)=T_{(\beta-\alpha)\gamma}(x,y)$ and $T(x,y)=T_{\beta \gamma}(x,y)$. It's clear that
    \[
        T=R+S,\ R(\beta,x)=T(\alpha,x),\ S(\beta,x)=T(\beta-\alpha,x)
    \]
    and that
    \[
        R(\beta,\beta),\ S(\beta,\beta),\ T(\alpha,\alpha),\ T(\beta-\alpha,\beta-\alpha)>0.
    \]
    Therefore, we can apply \cref{lem: estimate of norms} and obtain the desired inequality.
    In general, we have inequality
    \[
        G_{(1-\epsilon)\alpha \gamma,\beta}(x,x)\leq 5G_{\beta \gamma,\beta}(x,x),
    \]
    for $0<\epsilon<1$. Taking the limit $\epsilon\to 0$ supplies the desired inequality.
\end{proof}
By Poincar\'e duality, for any $h\in K$, we are allowed to define the dual norm of $G_{h^2}$ on $A^3$:
\[
    \|z\|_{h}^2=\sup_{0\neq x\in A^1}\frac{\deg(zx)^2}{G_{h^2,h}(x,x)}.
\]
Since the inner product $G_{h^2,h}$ is induced by the isomorphism between $A^1$ and $A^3\cong (A^1)^\vee $ given by the composition of multiplication by $h^2$ and $\theta_{T_{h^2},h}$, we see that
\[
    \|h^2x\|_h^2=G_{h^2,h}(\theta_{T_{h^2},h}(x),\theta_{T_{h^2},h}(x))=G_{h^2,h}(x,x),
\]
for any $x\in A^1$. We also have the following estimate of the norm $\|\cdot\|_{h}$.
\begin{lemma}\label{lem: estimate of external norm}
    Let $h,\alpha,\beta\in K$ such that $\alpha,\beta\preceq h$. Then we have
    \[
        \|\alpha\beta x\|_h^2\leq 25 G_{\alpha\beta,h}(x,x)
    \]
    for every $x\in A^1$.
\end{lemma}

\begin{proof}
    By applying \cref{lem: estimate of norms revisited} twice, we get
    \[
        G_{\alpha\beta,h}(y,y)\leq 5G_{h\beta,h}(y,y)\leq 25 G_{h^2,h}(y,y),
    \]
    and by \cref{lem: cauchy-schwarz},
    \[
        \deg(\alpha\beta x y)^2=T_{\alpha\beta}(x,y)^2\leq G_{\alpha\beta,h}(x,x)G_{\alpha\beta,h}(y,y)\leq 25 G_{\alpha\beta,h}(x,x)G_{h^2,h}(y,y).
    \]
    After dividing by $G_{h^2,h}(y,y)$ and taking the supremum over $y$, we obtain the desired inequality.
\end{proof}

Now we are ready to prove the main theorem using the Lefschetz decomposition and the above estimates.
\begin{proof}[Proof of \cref{thm: main}]
Let $R=\R[X_1,\ldots,X_n]$ and $P, H, V_i, A_i, B_i$ be as in the setup of \cref{thm: main}. Fix a graded ring homomorphism $R\to A$ that sends $K_R$ to $K$, and denote the images of $P,H,V_i,A_i,B_i$ by $p,h,v_i,\alpha_i,\beta_i$ respectively.

We first prove the theorem when $0\prec \alpha_i,\beta_i\preceq h$. According to the Lefschetz decomposition theorem, we can write $p$ uniquely as
\[
    p=h^2p_0+hp_1+p_2, \quad p_i\in A^i
\]
such that
\[
    h^3p_1=hp_2=0.
\]
By orthogonality of the Lefschetz decomposition, we have
\[
    \deg(p^2)=\deg(h^4p_0^2)+\deg(h^2p_1^2)+\deg(p_2^2)=p_0^2\deg(h^4)-G_{h^2,h}(p_1,p_1)+\deg(p_2^2).
\]
Meanwhile,
\[
    \|hp\|_h^2=\|h^2(hp_0+p_1)\|_h^2=G_{h^2,h}(hp_0+p_1,hp_0+p_1)=p_0^2\deg(h^4)+G_{h^2,h}(p_1,p_1).
\]
By adding the above equations and applying the Hodge-Riemann relation in degree $2$, we get
\[
    \deg(p^2)+\|hp\|_h^2=2p_0^2\deg(h^4)+\deg(p_2^2)\geq 2p_0^2\deg(h^4)=\frac{2\deg(ph^2)^2}{\deg(h^4)}.\tag{2}\label{eq: external inequality}
\]
Now apply the triangle inequality and \cref{lem: estimate of external norm} on the image of equality \eqref{eq: referencing}, we have
\[
    \|hp\|_h\leq \sum_{i=1}^r\|\alpha_i\beta_iv_i\|_h\leq\sum_{i=1}^r 5\sqrt {G_{\alpha_i\beta_i,h}(v_i,v_i)}=5\sum_{i=1}^r\sqrt{\frac{2\deg(\alpha_i\beta_ihv_i)^2}{\deg(\alpha_i\beta_ih^2)}-\deg(\alpha_i\beta_iv_i^2)}.
\]
By substituting the estimate above into inequality \eqref{eq: external inequality} we obtain the desired inequality.

In general, when $0\preceq\alpha_i,\beta_i\preceq h$ such that $\deg(\alpha_i\beta_ih^2)>0$, we can use a limit argument as follows: we can define perturbations of the chosen morphism $\varphi\colon R\to A$ as
\[
    \varphi_{\epsilon}\colon R\to A,\;X_i\mapsto\varphi(X_i)+\epsilon \ell
\]
where $\ell\in  K$ is any fixed ample class. Then for $0<\epsilon\ll 1$, $\varphi_\epsilon$ sends $\overline{K}_R\setminus0\to K$. In particular, if we denote the images of $P,H,V_i,A_i,B_i$ by $p_\epsilon,h_\epsilon,v_{i,\epsilon},\alpha_{i,\epsilon},\beta_{i,\epsilon}$, we have shown that
\[
   \deg(p_\epsilon^2)+25\left(\sum_{i=1}^r\sqrt{\frac{2\deg(\alpha_{i,\epsilon}\beta_{i,\epsilon}h_{\epsilon}v_{i,\epsilon})^2}{\deg(\alpha_{i,\epsilon}\beta_{i,\epsilon}h_{\epsilon}^2)}-\deg(\alpha_{i,\epsilon}\beta_{i,\epsilon}v_{i,\epsilon}^2)}\right)^2\geq \frac{2\deg(p_\epsilon h_\epsilon^2)^2}{\deg(h_\epsilon^4)},
\]
note that the assumption on the denominators ensures $A_i,B_i,H\neq 0$. Since both sides are continuous in $\epsilon$, after taking the limit of $\epsilon\to 0$, we obtain the desired inequality. Note that \[\deg(h^4)=\deg(\alpha_i\beta_ih^2)+\deg((h-\alpha_i)\beta_ih^2)+\deg(\alpha_i(h-\beta_i)h^2)+\deg((h-\alpha_i)(h-\beta_i)h^2)>0,\]
so the right-hand side has a well-defined limit.
\end{proof}

\begin{remark}
\begin{enumerate}[leftmargin=*]
    \item When $r=1$ and $A_1=B_1=H$, the theorem degenerates to
    \[\deg(h^2v_1^2)+25\left(\frac{2\deg(h^3v_1)^2}{\deg(h^4)}-\deg(h^2v_1^2)\right)\geq\frac{2\deg(h^3v_1)^2}{\deg(h^4)},\]
    which can be easily deduced from the AF-type inequality.
    \item 
    The estimate of the left-hand side of \eqref{eq: external inequality} is necessary. If $R\to A$ is not onto, the preimage of $hp$ under the Lefschetz isomorphism might not be contained in the image of $R$, and hence the value of $\|hp\|_h$ could be unknown. The current estimate is far from being sharp, but we'll see it is enough for the purpose of section \ref{sec: ex}.
    \item If we allow $A$ and $p$ to have higher degrees, one can still utilize the Lefschetz decomposition and repeatedly construct norms on $A^i$, and it's not hard to write down similar inequalities like \eqref{eq: external inequality}. However, obtaining analogous norm estimates like \Cref{lem: estimate of external norm} remains an open problem.
\end{enumerate}
\end{remark}
\section{Lorentzian polynomials and their apolar algebra}

We collect useful facts on Lorentzian polynomials in order to simplify the analysis in the next section.

Let $R=\R[X_1,\ldots,X_n]$ and $S=\R[x_1,\ldots,x_n]$ be two polynomial rings. Following the idea of Macaulay's \textit{inverse systems}, we define an action of $R$ on $S$ by
\[
    f(X_1,\ldots,X_n)\cdot Q\coloneqq  f(\frac{\partial}{\partial x_1},\ldots,\frac{\partial}{\partial x_n})Q,
\]
and we use $\ann(Q)$ to denote the ideal
\[
    \{f\in R\mid f\cdot Q=0\}. 
\]
When $0\neq Q\in S$ is homogeneous of degree $d\geq 2$, we define its \textit{apolar algebra} as the graded ring
\[
    A^\bullet_Q\coloneqq R/\ann(Q),
\]
which has Poincar\'e duality with $\deg\colon A^d_Q\to \R$ given by $\deg(f)=f\cdot Q$.

Let $(A^\bullet,\deg)$ be a Poincar\'e duality algebra and elements $\chi_1,\ldots,\chi_n$ be elements in $A^1$. We define the \textit{volume polynomial} with respect to $\chi_1,\ldots,\chi_n$ by
\begin{align*}
    Q(x_1,\ldots,x_n)=&\frac{1}{d!}\deg\big((x_1\chi_1+\cdots+x_n\chi_n)^d\big)\\
    =&\sum_{\substack{
    k_1+\cdots +k_n=d\\k_1,\ldots ,k_n\geq 0    }}\deg(\chi_1^{k_1}\cdots\chi_n^{k_n})\frac{x_1^{k_1}\cdots x_n^{k_n}}{k_1!\cdots k_n!}.
\end{align*}
We have the following correspondence.
\begin{theorem}[see for example, {\cite[Theorem 21.6 and Exercise 21.7]{Eis95}}] The map
\[
    \big\{0\neq Q\in S^d\big\}\xrightarrow{Q\mapsto A^\bullet_Q}\big\{\substack{\text{graded quotient rings of}\ R\ \text{with a degree map in degree}\ d\\ \text{satisfying Poincar\'e duality}}\big\},
\]
is bijective, with the inverse given by taking the volume polynomial with respect to the images of $X_1,\ldots,X_n$.
\end{theorem}

A homogeneous polynomial of degree $d$ is \textit{strictly Lorentzian} if it has positive coefficients and the Hessian matrices of its $(d-2)$-th order partial derivatives have the Lorentzian signature; a polynomial is \textit{Lorentzian} if it's the limit of strictly Lorentzian polynomials. For a homogeneous polynomial $Q\neq 0$, we define $K_Q\coloneqq \{a_1X_1+\cdots+a_nX_n\mid a_1,\ldots,a_n>0\}\subseteq A^1_Q$. We say a Poincar\'e duality algebra satisfies \textit{mixed $(\mathrm{HR}^{\leq m})$} with respect to a cone $K$ if the mixed Hodge-Riemann relations are satisfied for $0\leq k\leq m$. We have the following characterization of Lorentzian polynomials:

\begin{theorem}[{\cite[Theorem 3.2]{Eur21}}]\label{thm: lorentzian and HR}
A nonzero homogeneous polynomial $Q\in S$  is Lorentzian if and only if the triple $(A^\bullet_Q,\deg,K_Q)$ satisfies mixed $(\mathrm{HR}^{\leq 1})$.
\end{theorem}

Let $d\Delta_n=\{k=(k_1,\ldots,k_n)\in \mathbb Z_{\geq 0}^n\mid \sum_{i=1}^n k_i=d\}$. For $Q\in S^d$, the support of $Q$, denoted $\supp Q$, is the set of exponent vectors where the corresponding coefficient of $Q$ is nonzero. A subset $\mathcal S$ of $d\Delta_n$ is \textit{M-convex} if for any $k,k'\in \mathcal S$ and for any $i$ such that $k_i<k'_i$, there exists $j$ such that $k_j>k'_j$ and both
\[
    k+e_i-e_j\text{ and } k'-e_i+e_j
\]
are contained in $\mathcal S$. We also need the following characterization of Lorentzian polynomials.

\begin{theorem}[{\cite[Theorem 2.25]{BH20}}]\label{thm: Lorentzian-criterion}
A homogeneous polynomial $Q\in S^d$ is Lorentzian if and only if
\begin{enumerate}
    \item  Q has nonnegative coefficients,
    \item $\supp(Q)\subseteq d\Delta_n$ is M-convex, and
    \item For any $k_1+\cdots +k_n=d-2,\ k_1,\ldots,k_n\geq 0$, the Hessian matrix of \[\left(\frac{\partial}{\partial x_1}\right)^{k_1}\cdots\left(\frac{\partial}{\partial x_n}\right)^{k_n}Q \] has at most one positive eigenvalue.
\end{enumerate}
\end{theorem}

We say a polynomial $Q\in S$ is \textit{concise} if $R^1\cong A^1_Q$, i.e., $f\cdot Q\neq 0$ for any nonzero linear form $f\in R$. The following variation of \cite[Theorem 2.10]{BH20} would be useful for proving the strict Lorentzian property.

\begin{proposition}\label{prop: strictification}
   Let $(A^\bullet,\deg,K)$ be a Poincar\'e duality algebra satisfying \textit{mixed $(\mathrm{HR}^{\leq 1})$} with respect to $K$.
   \begin{enumerate}
   \item If $\gamma_1,\ldots,\gamma_m\in\overline{K}$, then their volume polynomial is Lorentzian;
   \item If $\gamma_1,\ldots,\gamma_m\in K$ are linearly independent, then their volume polynomial is strictly Lorentzian.
   \end{enumerate}
   In particular, for any rank $m$ real $n\times m$ matrix $M$ with positive entries and any concise Lorentzian polynomial $Q$ in $n$ variables, $Q(My)$ is a strict Lorentzian in $m$ variables.
\end{proposition}
\begin{proof}
    Consider the case when $\gamma_1,\ldots,\gamma_m\in K$. Let $F=F(y_1,\ldots,y_m)$ be the volume polynomial of $\gamma_1,\ldots,\gamma_m$. Mixed Hodge-Riemann relation in degree 0 implies
\[
    \deg\!\left(\prod_{i=1}^m\gamma_i^{k_i}\right)>0
\]
for any $k_1+\cdots+k_m=d,\ k_1,\ldots,k_m\geq 0$. In other words, $F$ has positive coefficients.

In degree 1, the Hodge-Riemann relation implies that the symmetric form
\[
    A^1\times A^1 \to \R, \ \ (\alpha,\beta)\mapsto \deg\left(\prod_{i=1}^m\gamma_i^{k_i}\cdot \alpha\beta\right) 
\]
has the Lorentzian signature for any $k_1+\cdots +k_m=d-2,\ k_1,\ldots,k_m\geq 0$. Cauchy's interlacing formula, plus the fact that this form is positive on $\gamma_i$, implies that the restriction of this form to the subspace $\operatorname{span}\{ \gamma_i\}_i$ also has the Lorentzian signature. By pulling back this form to $\R^m$ via $e_i\mapsto\gamma_i$, we obtain the Hessian matrix of $\left(\frac{\partial}{\partial y_1}\right)^{k_1}\cdots\left(\frac{\partial}{\partial y_m}\right)^{k_m}F$, so this matrix has exactly one positive eigenvalue. This proves that $F$ is Lorentzian, by \ref{thm: Lorentzian-criterion}.
When $\gamma_1,\cdots,\gamma_m\in \overline{K}$, a limiting argument shows that $F$ is still Lorentzian; when $\gamma_1,\ldots,\gamma_m\in K$ are linearly independent, one further sees that the Hessian matrices have the Lorentzian signature and therefore $F$ is strictly Lorentzian.
    Given a concise Lorentzian polynomial $Q$, consider its apolar algebra with the standard generators $\chi_1,\ldots,\chi_n$. Let $\gamma_1,\ldots,\gamma_m$ be elements in $A^1$ defined by $\gamma=M^T\chi$. One sees that 
    \[
    \deg(f(\gamma))=\deg(f(M^T\chi))=f(M^TX)\cdot Q(x)=f(Y)\cdot Q(My)
    \]
    for any homogeneous polynomial $f$ of degree $d$, where $X_i$ and $Y_j$ act as partial differentials. Since $M$ has positive entries and that $M$ is of rank $m$, $\gamma$'s are linearly independent elements in $K$ and $Q(My)$ is precisely their volume polynomial. Therefore, $Q(My)$ is strictly Lorentzian by the previous claim.

\end{proof}

\begin{remark}
    The proposition above also shows that $A^\bullet_{Q(My)}$ satisfies mixed Hodge-Riemann relations in degree less than or equal to $1$ for classes in $\overline K_{Q(My)}\setminus 0$, by \cref{thm: lorentzian and HR}. This implication cannot extend to higher Hodge-Riemann relations. See \cite[Example 6.17]{NW26}, for example.
\end{remark}

\section{Analysis of examples}\label{sec: ex} 
In this section, we verify the various claims made in \cref{ex: ternary}. Let $a,b,r,s, F$ and $A, \chi,\gamma,\zeta$ be as in the setup of \cref{ex: ternary}.
\subsection{Failure of the inequality} Notice that $F$ is symmetric in $y,z$. Let $\tilde{p}=\frac{p}{750}=\chi^2-\chi(\gamma+\zeta)-\frac{14}{375}(\gamma^2+\zeta^2)+\frac{49}{50}\gamma\zeta$. When $\epsilon=0$, we have
\begin{align*}
    \deg(\tilde p \chi^2)=&-\frac{7}{250}\\
    \deg(\tilde p \chi\gamma)=\deg(\tilde p\chi\zeta)=&\frac{7}{500}\\
    \deg(\tilde p \gamma^2)=\deg(\tilde p \zeta^2)=&-\frac{17}{300}\\
    \deg(\tilde p \gamma\zeta)=&
\frac{16}{375}\end{align*}
and therefore
\[
    \deg(\tilde p^2)=(-\frac{7}{250})-2\cdot \frac{7}{500}-2\cdot\frac{14}{375}\cdot(-\frac{17}{300})+\frac{49}{50}\cdot \frac{16}{375}=-\frac{56}{5625}
\]
and
\[
    \deg(p^2)=750^2\cdot(-\frac{56}{5625})=-5600.
\]
It now suffices to compute $G_{\alpha_i\beta_i,h}(v_i,v_i)$, where these classes are given by the right-hand side of the decomposition \eqref{eq: referencing'}. We note that since $F$ has positive coefficients, we can readily verify the denominators appearing in $G_{\alpha_i\beta_i,h}(v_i,v_i)$ are positive. To simplify notations,  for a quadratic $Q$, we denote the Hessian matrix of $Q\cdot F$ by $H_Q$, then we have
\[
H_{X^2}=
\begin{pmatrix}
a&b&b\\
b&1&r\\
b&r&1
\end{pmatrix},\quad
H_{Y^2}=
\begin{pmatrix}
1&s&1\\
s&s^2&s\\
1&s&1
\end{pmatrix},\quad
H_{Z^2}=
\begin{pmatrix}
1&1&s\\
1&1&s\\
s&s&s^2
\end{pmatrix},
\]
and
\[
H_{XY}=
\begin{pmatrix}
b&1&r\\
1&s&1\\
r&1&1
\end{pmatrix},\quad
H_{XZ}=
\begin{pmatrix}
b&r&1\\
r&1&1\\
1&1&s
\end{pmatrix},\quad
H_{YZ}=
\begin{pmatrix}
r&1&1\\
1&s&1\\
1&1&s
\end{pmatrix}.
\]
Direct computation gives
\begin{gather*}
    H_{X^2}\begin{bmatrix}
        30\\-17\\-17
    \end{bmatrix}=\vec 0, \ \ H_{XY}\begin{bmatrix}
       10 \\ -2 \\ -9 
    \end{bmatrix}=\vec 0,\ \ 
    H_{XZ}\begin{bmatrix}
       10 \\ -9 \\ -2
    \end{bmatrix}=\vec 0,\\
    H_{Y^2}\begin{bmatrix}
       1 \\ 0 \\ -1
    \end{bmatrix}=H_{Y^2}\begin{bmatrix}
       0 \\ 2 \\ -1
    \end{bmatrix}=\vec 0,\ \ 
    H_{Z^2}\begin{bmatrix}
       1 \\ -1 \\ 0
    \end{bmatrix}=H_{Z^2}\begin{bmatrix}
       0 \\ -1 \\ 2
    \end{bmatrix}=\vec 0,
\end{gather*}

where the first equality implies $\deg(\chi^2(30\chi-17\gamma-17\zeta)w)=0$ for any $w\in \mathrm{span}\{\chi,\gamma,\zeta\}$, and in particular,
\[
    G_{\chi^2,h}(30\chi-17\gamma-17\zeta,30\chi-17\gamma-17\zeta)=0.
\]
The same argument works for other equalities, and we have that
\[
    \deg(p^2)+25\left(\sum_{i=1}^r\sqrt{G_{\alpha_i\beta_i,h}(v_i,v_i)}\right)^2=\deg(p^2)=-5600<0\leq \frac{2\deg(ph^2)^2}{\deg(h^4)}
\]
when $\epsilon=0$. Since the left-hand side is continuous in $\epsilon$, the denominators stay positive and the inequality stays true for $0<\epsilon\ll 1$.

\subsection{Verification of the Lorentzian property}
Now we show that $F$ is Lorentzian and that $F_\epsilon$ is strictly Lorentzian for $0<\epsilon\ll 1$. Since $F_{\epsilon}=F(Mx)$ where $M=I+\epsilon \begin{bsmallmatrix}
    1&1&1\\1&1&1\\1&1&1
\end{bsmallmatrix}$ is invertible with positive entries for $0<\epsilon \ll 1$, by \cref{prop: strictification} it suffices to prove that $F$ is concise and Lorentzian. Since $F$ has full support, for Lorentzian property, it remains to check that the various matrices $H_Q$ have at most one positive eigenvalue, according to \cref{thm: Lorentzian-criterion}.

For $Q=X^2,XY$ or $XZ$, we've seen that $H_Q$ has one null vector, given by the factors on the right-hand side of \eqref{eq: referencing'}, so it suffices to show that $H_Q$ is not positive semi-definite. This can be done by evaluating principal minors, and we have
\[
    \det(H_{X^2,\{2,3\}})=1-r^2<0,\ \ \det(H_{XY,\{2,3\}})=\det(H_{XZ,\{2,3\}})=s-1<0,
\]
where $H_{Q,I}$ denotes the principal submatrix of $H_Q$ by taking columns and rows in $I$. For $Q=Y^2$ or $Z^2$, there's nothing to check since the null space of $H_Q$ has dimension $2$. For $Q=YZ$, we have
\[
    H_{YZ}\begin{bmatrix}
        0\\1\\-1
    \end{bmatrix}=-\frac{1}{2}\begin{bmatrix}
        0\\1\\-1
    \end{bmatrix}
\]
and
\[
    \det(H_{YZ})=\frac{7}{40}>0,
\]
so $H_{YZ}$ has the Lorentzian signature.
    Now it remains to check conciseness. For any linear form $L$ such that $L\cdot F=0$, we must also have $L\cdot(YZ\cdot F)=0$. But we just found that $H_{YZ}$ is non-singular, therefore $L=0$.

\begin{remark}
    The first breakthrough of this project was a similar example given by the Danus system to illustrate that Lorentzian polynomials might fail to arise from the K\"ahler package. The original argument was to show that the element $p$ is truly primitive with respect to $h$ and therefore $\deg(p^2)$ should be nonnegative. The author later prompted the system to transform the argument into an inequality. The example presented here is also generated based on the same idea but is easier to present.
\end{remark}

\subsection{Comparison of types of polynomials}

We investigate different sets of polynomials. Let $L^d_n\subseteq S^d$ ($\mathring{L}^d_n$) denote the set of (strictly) Lorentzian polynomials in $n$ variables of degree $d$. According to \cref{thm: lorentzian and HR}, these are precisely the polynomials whose apolar algebras have mixed ($\text{HR}^{\leq 1}$) when nonzero. We follow the idea of \cite{MMS25} and make the following definition.

\begin{definition}
    A homogeneous polynomial $F\in S^d$ is said to be $k$-Lorentzian if $F=0$ or the triple $(A^\bullet_F,\deg,K_F)$ satisfies mixed ($\text{HR}^{\leq k}$). The set of all $k$-Lorentzian polynomials in $n$ variables in degree $d$ is denoted $L(k)^d_n$. Note that $L(1)^d_n=L^d_n$.
\end{definition}

\begin{definition}
    Let $\Bbbk$ be an algebraically closed field.
    \begin{enumerate}
        \item A homogeneous polynomial $F\in S^d$ of degree $d$ is called a \textit{realizable volume polynomial over $\Bbbk$} if there exists a $d$-dimensional projective variety $Y$ over $\Bbbk$ and a collection of semiample divisors $D_1,\ldots,D_n$ on $Y$ such that
    \[
        F=\frac{q}{d!}\int_Y\bigg(\sum_{i=1}^nx_iD_i\bigg)^d
    \]
    for some $q\in\Q_{\geq 0}$. The set of all such polynomials is denoted $\Bbb V^d_n(\Q,\Bbbk)$
        \item A homogeneous polynomial $F\in S^d$ is a \textit{volume polynomial over $\Bbbk$} if it is a limit of realizable volume polynomials over $\Bbbk$. The set of all such polynomials is denoted $\Bbb V^d_n(\R,\Bbbk)$.
    \end{enumerate}
\end{definition}

We introduce new families of polynomials as follows. Note that there's no need to invoke (mixed) hard Lefschetz, see for example \cite[Lemma 2.10]{MMS25}.

\begin{definition}
    \begin{enumerate}[leftmargin=*]
        \item A homogeneous polynomial $F\in S^d$ is a \textit{strict volume polynomial with respect to mixed $(\text{HR}^{\leq k})$} if there exists a Poincar\'e duality algebra $(A^\bullet,\deg,K)$ satisfying mixed $(\text{HR}^{\leq k})$, and there exists $\ell_1,\ldots,\ell_n\in K$, linearly independent, such that $F$ is the volume polynomial with respect to $\ell_1,\ldots,\ell_n$. The set of all such polynomials is denoted $\mathring{\Bbb V}^d_{\text{HR}^{\leq k},n}$
        \item A homogeneous polynomial $F\in S^d$ is a \textit{volume polynomial with respect to mixed $(\text{HR}^{\leq k})$} if there exists a Poincar\'e duality algebra $(A^\bullet,\deg,K)$ satisfying mixed $(\text{HR}^{\leq k})$, and there exists $\ell_1,\ldots,\ell_n\in \overline K$ such that $F$ is the volume polynomial with respect to $\ell_1,\ldots,\ell_n$. The set of all such polynomials is denoted $\Bbb V^d_{\text{HR}^{\leq k},n}$. We also let $\Bbb {V'}^{d}_{\text{HR}^{\leq k},n}$ denote the subset of $\Bbb {V}^{d}_{\text{HR}^{\leq k},n}$ such that the algebra can be chosen so that $n\leq\dim A^1$. Note that this automatically includes all concise polynomials in $\Bbb {V}^{d}_{\text{HR}^{\leq k},n}$.
    \end{enumerate}
\end{definition}

\begin{proposition}\label{prop: volume and lorentzian}
    A (strict) volume polynomial with respect to mixed $(\text{HR}^{\leq 1})$ is the same as a (strictly) Lorentzian polynomial.
\end{proposition}

\begin{proof}
    One direction is supplied by \cref{prop: strictification}.

    For the other direction, \Cref{thm: lorentzian and HR} shows that $\Bbb V^d_{\text{HR}^{\leq 1},n}\supseteq L^d_n$. Given a strictly Lorentzian polynomial $F(x)$, the polynomial $G(x)=F((I+\epsilon E)^{-1}x)$ is also strictly Lorentzian for $0<\epsilon\ll1$, where $E$ is the all-ones matrix. It's clear that in the apolar algebra of $G$, the volume polynomial with respect to $X_1+\epsilon(X_1+\cdots+X_n),X_2+\epsilon(X_1+\cdots+X_n),\ldots,X_n+\epsilon(X_1+\cdots+X_n)$ is exactly $F$, which shows that $F$ is a strict volume polynomial with respect to mixed $(\text{HR}^{\leq 1})$.
\end{proof}

We have the following comparison of different notions.
\begin{proposition}
    We have the following inclusions for any $1\leq k\leq \lfloor\frac{d}{2}\rfloor$.
    \[\begin{tikzcd}
	 & & L^d_n & \\
	 \Bbb V^d_n(\R,\Bbb C) & \overline{\mathring{\Bbb V}^d_{\text{HR}^{\leq k},n}} & {\Bbb V}^d_{\text{HR}^{\leq k},n} & L(k)^d_n \\
    \Bbb V^d_n(\Q,\Bbb C) & {\Bbb V'}^d_{\text{HR}^{\leq k},n} & & \\
	& \mathring{\Bbb V}^d_{\text{HR}^{\leq k},n}
	\arrow[hook, from=2-1, to=2-2]
	\arrow[hook, from=2-2, to=1-3]
	\arrow[hook, from=2-3, to=1-3]
	\arrow[hook', from=2-4, to=2-3]
	\arrow[hook, from=3-1, to=2-1]
    \arrow[hook, from=3-1, to=3-2]
	\arrow[hook, from=3-2, to=2-2]
	\arrow[hook, from=3-2, to=2-3]
	\arrow[hook, from=4-2, to=3-2]
\end{tikzcd}\]
All inclusions here, except for the ones starting from ${\Bbb V'}^d_{\text{HR}^{\leq k},n}$, are strict in general.
\end{proposition}

\begin{proof}
    To check the inclusion from the left column to the middle column, given a realizable volume polynomial, after blowing up the variety $Y$ enough times, one can take $\bigoplus_kH^{k,k}(Y)\cap H^\bullet(Y,\R)$ which has the K\"ahler package and a sufficiently large degree $1$ component. Other inclusions directly follow from the definition. The vertical inclusion on the left is strict simply because a realizable volume polynomial cannot have irrational coefficients.

    According to \cite[Footnote 15]{BH20}, there exists a concise cubic Lorentzian polynomial that fails to satisfy the \textit{reverse Khovanskii-Teissier inequality}, hence fails to be a volume polynomial over any field. On the other hand, \cref{prop: volume and lorentzian} shows that it is a volume polynomial with respect to mixed $(\text{HR}^{\leq 1})$, so the inclusions from the left column to the middle column are strict when $k=1$.
    
    The existence of concise non-strict Lorentzian polynomial, for example $xy$, shows the strictness of the vertical inclusion at the bottom when $k=1$.

    As we've seen in \cref{ex: ternary}, a strictly Lorentzian polynomial (and its small perturbations) might not be realized as the volume polynomial in any ring with mixed $(\text{HR}^{\leq 2})$, which shows the strictness of the inclusions $\overline{\mathring{\Bbb V}^d_{\text{HR}^{\leq k},n}} \subseteq L^d_n$ and ${\Bbb V}^d_{\text{HR}^{\leq k},n} \subseteq L^d_n$ when $k\geq 2$.

    \cite[Example 6.17]{NW26} shows that there exist binary quartic realizable volume polynomials whose apolar algebra doesn't satisfy mixed $(\text{HR}^{\leq 2})$, which means the inclusions $L(k)^d_n\subseteq {\Bbb V}^d_{\text{HR}^{\leq k},n}$ and $L(k)^d_n\subseteq L^d_n$ are strict for $k\geq 2$.
\end{proof}

\begin{remark}
\begin{enumerate}[leftmargin=*]
    \item The introduction of ${\Bbb V'}^d_{\text{HR}^{\leq k},n}$ concerns whether one can enlarge an algebra with the K\"ahler package without changing its degree. When $k=\lfloor\frac{d}{2}\rfloor$, one can consider the algebras
    \[
        B\coloneqq A\otimes \R[x]/x^2,\quad C\coloneqq B/\ann(h+x),
    \]
    where $h\in K$. By \cite[Proposition 7.7]{AHK18} and \cite{Cat08}, $B$ has the K\"ahler package with respect to $(K+\R_{>0}x)$. Therefore, by the descent lemma \cref{prop: descent}, $C$ also has the K\"ahler package with respect to the image of $(K+\R_{>0}x)$. When $d\geq 3$, we have
    \[
        \dim C^1=\dim A^1+1,\qquad \deg_C(\bar \alpha_1\cdots\bar \alpha_d)=\deg_A(\alpha_1\cdots\alpha_d),
    \]
    for elements $\alpha_1,\ldots\alpha_d\in A^1$. This shows that ${\Bbb V}^d_{\text{HR}^{\leq k},n}={\Bbb V'}^d_{\text{HR}^{\leq k},n}$ and we can compare $\overline{\mathring{\Bbb V}^d_{\text{HR}^{\leq k},n}}$ and ${\Bbb V}^d_{\text{HR}^{\leq k},n}$ directly. For general $k$, it remains open whether this ``blowup'' construction still works.
    \item When $k=1$, \Cref{prop: volume and lorentzian} shows that $\overline{\mathring{\Bbb V}^d_{\text{HR}^{\leq k},n}}$ and ${\Bbb V}^d_{\text{HR}^{\leq k},n}$ are equal. Whether this remains true when $k>1$ remains open.
\end{enumerate}
\end{remark}

\bibliographystyle{amsalpha}
\bibliography{references}
\end{document}